\documentclass[12pt]{amsart}
\usepackage{tikz-cd}
\usepackage[T1]{fontenc}
\usepackage{lmodern}
\usepackage{microtype}
\usepackage[a4paper,margin=31mm]{geometry}
\usepackage{mathtools}
\usepackage{amssymb}
\usepackage{xcolor}
\usepackage{enumitem}
\usepackage{hyperref}
\hypersetup{colorlinks=true,linkcolor=blue!45!black,
citecolor=blue!45!black,urlcolor=blue!55!black, pdftitle={Disjoint Divisors and Splitting Rigidity under Nonnegative Orthogonal Ricci Curvature}}
\setlist[itemize]{leftmargin=2em,itemsep=2pt,topsep=4pt}
\setlist[enumerate,1]{label=\textup{(\roman*)},leftmargin=2em,itemsep=2pt,topsep=4pt}
\setlist[enumerate,2]{label=\textup{(\alph*)},leftmargin=2em,itemsep=2pt,topsep=4pt}
\numberwithin{equation}{section}
\theoremstyle{plain}
\newtheorem{theorem}{Theorem}[section]
\newtheorem{proposition}[theorem]{Proposition}
\newtheorem{lemma}[theorem]{Lemma}
\newtheorem{corollary}[theorem]{Corollary}
\newtheorem{definition}{Definition}[section]
\newtheorem{example}[theorem]{Example}
\theoremstyle{remark}

\newcommand{\Pj}{\mathbb P}
\newcommand{\Ric}{\operatorname{Ric}}
\newcommand{\rd}{\operatorname{rd}}
\newcommand{\Supp}{\operatorname{Supp}}
\newcommand{\Isomhol}{\operatorname{Isom}_{\mathrm{hol}}}
\newcommand{\II}{\mathrm{II}}
\newcommand{\cO}{\mathcal O}
\newcommand{\cV}{\mathcal V}

\title[Disjoint Divisors and Splitting Rigidity under Nonnegative Orthogonal Ricci Curvature]{Disjoint Divisors and Splitting Rigidity under Nonnegative Orthogonal Ricci Curvature}
\author{Xingluan Wang}
\date{}
\thanks{E-mail: wxlmath@mail.ustc.edu.cn. School of Mathematical Sciences, University of Science and Technology of China, Hefei 230026, P. R. China.}


\begin{document}
\raggedbottom
\begin{abstract}
Let $X^n$ be a compact K\"ahler manifold with nonnegative orthogonal Ricci curvature.
We prove that if $X$ admits two nonzero effective divisors with disjoint supports, then its universal cover splits.
As applications, we study the MRC fibration of $X$, and show that when the dimension of a general fiber is $n-1$, the universal cover of $X$ splits as a product of a simply connected Riemann surface and a compact rationally connected manifold.
When $X$ is rationally connected, it splits off a $\Pj^1$-factor if and only if it contains two disjoint nonzero effective divisors.
Together with the results of Zhang and Zhang, these conclusions give a unified description of the universal cover for all possible dimensions of the general MRC fiber.
We also establish metric existence criteria for certain projective bundles.
\end{abstract}
\maketitle
\pagestyle{plain}

\section{Introduction}\label{intro}
\subsection{Background and main results}
For a unit vector $v\in T^{1,0}X$ on a K\"ahler manifold $(X,g)$, the
\emph{orthogonal Ricci curvature} is defined by
\[
\Ric_g^\perp(v)
=\Ric_g(v,\bar v)-R_g(v,\bar v,v,\bar v).
\]
Ni and Zheng systematically studied this curvature, establishing comparison and vanishing theorems and proving that compact Kähler manifolds with positive orthogonal Ricci curvature are projective \cite{NiZheng}. Further developments can be found in \cite{Ni2021,NiWangZheng,ZhangZhang}.

The purpose of this paper is to establish splitting results for compact
K\"ahler manifolds with nonnegative orthogonal Ricci curvature. We use three related notions of splitting structures in this paper.

\begin{definition}\label{ks}
The universal cover of a K\"ahler manifold $(X,g)$ is said to
\emph{split holomorphically and isometrically} if it admits a product
decomposition
\[
(\widetilde X,\widetilde g)
\simeq (X_1,g_1)\times(X_2,g_2)
\]
by complete simply connected K\"ahler manifolds.

A five-tuple $(X,Y,F,f,\rho)$ is a \emph{locally constant fibration} if
$X,Y,F$ are compact complex manifolds, $f\colon X\to Y$ is a surjective
holomorphic map with connected fibers, and there is a representation
$\rho\colon\pi_1(Y)\to\operatorname{Aut}(F)$ for which
\[
X\simeq(\widetilde Y\times F)/\pi_1(Y)
\]
biholomorphically over $Y$, where
$\gamma\cdot(y,p)=(\gamma y,\rho(\gamma)p)$.

If, moreover, $X,Y,F$ carry K\"ahler metrics $g,g_Y,g_F$ such that
$\rho(\pi_1(Y))\subseteq\Isomhol(F,g_F)$ and the above identification
is a holomorphic isometry
\[
(X,g)\simeq
\bigl((\widetilde Y,\widetilde g_Y)\times(F,g_F)\bigr)/\pi_1(Y),
\]
we call the locally constant fibration $(X,Y,F,f,\rho)$ a \emph{K\"ahler suspension}.
\end{definition}

In this paper, the fiber $F$ is rationally connected and hence simply connected.
A K\"ahler  suspension identification yields a holomorphic isometric splitting.

Splitting structures of these types are known under several
nonnegative curvature conditions. Cheeger and Gromoll established
Riemannian isometric splitting under nonnegative Ricci curvature. In
the compact K\"ahler setting their theorem yields a holomorphic
isometric splitting of the universal cover \cite{CheegerGromoll}.
Mok proved a more precise holomorphic isometric decomposition under
nonnegative holomorphic bisectional curvature \cite{Mok}, while Ni
and Tam obtained related splitting theorems in the complete noncompact
setting \cite{NiTam}.
Matsumura showed that compact K\"ahler manifolds with
semipositive holomorphic sectional curvature admit a K\"ahler
suspension with rationally connected fiber and flat base
\cite{Matsumura2026}.
Chu, Lee, and Zhu established universal-cover splitting results
under nonnegative mixed curvature \cite{ChuLeeZhu}.

Recently, Zhang and Zhang proved a splitting structure theorem in part of the
$\Ric^\perp \geq 0$ case \cite{ZhangZhang}. To describe their result, recall that the maximal rationally connected
(MRC) fibration is an almost holomorphic dominant meromorphic map whose
general fibers are rationally connected and whose base is not uniruled.
The \emph{rational dimension} $\rd(X)$ is the dimension of a general
MRC fiber. They proved that: 
if $(X^n,g)$ is a compact K\"ahler manifold with
$\Ric_g^\perp\geq0$ and $\rd(X)<n-1$, 
then $X$ admits a K\"ahler suspension
$(X,Y,F,f,\rho)$ where $F$ is rationally connected and $g_Y$ is flat \cite[Theorem~4.13]{ZhangZhang}.

Their Bochner argument does not apply when
$\rd(X)=n-1$. Our approach is motivated by the geometry of divisors. Ni, Wang, and Zheng proved that any two irreducible divisors on a compact Kähler manifold with \(\operatorname{Ric}^{\perp}>0\) must intersect \cite{NiWangZheng}. When the curvature is nonnegative, disjoint divisors may exist. Our first result
shows that their existence forces a global splitting.

\begin{theorem}\label{main}
Let $(X^n,g)$ be a compact K\"ahler manifold with $\Ric_g^\perp\geq0$, and let $D_1,D_2$ be nonzero effective divisors with disjoint supports.
Then the following statements hold.
\begin{enumerate}
\item The reduced divisors $(D_i)_{\mathrm{red}}$ are smooth and tangent to a unique parallel holomorphic hyperplane distribution $\cV$.
\item The universal cover admits a holomorphic isometric splitting
\[
\begin{gathered}
(\widetilde X,\widetilde g)
\simeq
(\Sigma,h_\Sigma)\times(N,g_N),\quad
\dim_{\mathbb C}\Sigma=1,
\quad
\widetilde{\cV}=T^{1,0}N,
\end{gathered}
\]
where both factors are complete and simply connected.
\item Each component of $(D_i)_{\mathrm{red}}$ is totally geodesic, and its normal line bundle has flat Chern connection with respect to the quotient metric.
\end{enumerate}
\end{theorem}

The existence of disjoint divisors is sufficient for this splitting, but the converse fails in general; see Example~\ref{count1}.
Theorem~\ref{main} also extends the divisor intersection theorem of Ni, Wang, and Zheng \cite[Theorem~1.6]{NiWangZheng} to the quasi-positive case.

\begin{corollary}\label{quasipositive}
Let $(X^n,g)$ be a compact K\"ahler manifold with $\Ric_g^\perp\geq0$.
Suppose that there exists a point $p\in X$ such that
\[
\Ric_g^\perp(v)>0
\qquad\text{for every unit vector }v\in T_p^{1,0}X.
\]
Then any two nonzero effective divisors on $X$ have intersecting supports.
\end{corollary}

Theorem~\ref{main} applies to the case $\rd(X)=n-1$:
the MRC map has a curve as its base, and two distinct smooth general
fibers provide disjoint divisors. When \(\rd(X)=n\), the existence of two disjoint nonzero effective divisors is equivalent to a holomorphic isometric splitting off a \(\mathbb P^1\)-factor.

\begin{theorem}\label{mrcsplit}
Let $(X^n,g)$ be a compact K\"ahler manifold with $\Ric_g^\perp\geq0$.
\begin{enumerate}
\item If $\rd(X)=n-1$, then the MRC fibration is represented by a holomorphic submersion $f\colon X\to C$ onto a smooth compact curve of genus at least one, and $X$ admits a K\"ahler suspension $(X,C,F,f,\rho)$ with rationally connected fiber $F$.
\item If $\rd(X)=n$, then the following conditions are equivalent:
\begin{enumerate}
\item $X$ contains two nonzero effective divisors with disjoint supports;
\item there exists a nonconstant holomorphic map $X\to\Pj^1$;
\item $(X,g)$ is holomorphically isometric to $(\Pj^1,h)\times(F,g_F)$, where $h,g_F$ are K\"ahler metrics and $F$ is compact and rationally connected.
\end{enumerate}
\end{enumerate}
\end{theorem}
When $\rd(X)<n-1$, the induced base metric is flat, whereas for $\rd(X)=n-1$ the base curve need not admit any flat K\"ahler metric; see Example~\ref{nonflat}.
Together with the theorem of Zhang and Zhang quoted above, this yields a uniform
description of the universal cover for all values of $\rd(X)$.

\begin{theorem}\label{universalstructure}
Let $(X,g)$ be a compact K\"ahler manifold with $\Ric_g^\perp\geq0$.
Then its universal cover admits a holomorphic isometric splitting
\[
(\widetilde X,\widetilde g)
\simeq(M,g_M)\times(F,g_F),
\]
where $F$ is a compact rationally connected K\"ahler manifold, and $(M,g_M)$ is either a point, a flat complex Euclidean space, or a complete simply connected noncompact Riemann surface.
\end{theorem}

In particular, our result rules out the higher-dimensional noncompact factor in \cite[Theorem~1.3(2)]{GuZhang}.
\subsection{Further applications}
We also establish existence criteria of metrics with nonnegative orthogonal Ricci curvature on some projective bundles.

\begin{proposition}
Let $M$ be a compact K\"ahler manifold, let $L$ be a holomorphic line bundle on $M$, and set $X=\Pj(\cO_M\oplus L)$.
Then the following conditions are equivalent:
\begin{enumerate}
\item $X$ admits a K\"ahler metric with $\Ric^\perp\geq0$;
\item $M$ admits a K\"ahler metric with $\Ric^\perp\geq0$, and $L$ is unitary flat.
\end{enumerate}
\end{proposition}

Apostolov, Calderbank, Gauduchon, and T{\o}nnesen-Friedman
\cite{ACGTF2011} proved that, for a compact Riemann surface
$C$ of genus at least two, the projective bundle $\mathbb P(E)$
admits a constant scalar curvature K\"ahler metric if and only
if $E$ is polystable.
We obtain an analogous existence criterion for K\"ahler metrics
with nonnegative orthogonal Ricci curvature.

\begin{theorem}
\label{projectivebundlescurves}
Let $C$ be a compact Riemann surface, and let
$E\to C$ be a holomorphic vector bundle of rank $r\geq2$.
Set $X=\mathbb P(E)$.
Then the following conditions are equivalent:
\begin{enumerate}
\item $E$ is polystable.
\item $X$ admits a K\"ahler metric $g$ with nonnegative orthogonal
Ricci curvature.
\end{enumerate}
\end{theorem}

\subsection{Organization of the paper}

Section~\ref{divisors} establishes the self-intersection inequality for effective divisors (Proposition~\ref{Dgeq0}) and the parallel-kernel criterion (Lemma~\ref{parallelkernel}), and combines them to prove Theorem~\ref{main}.
Section~\ref{fibrations} develops the passage from parallel distributions to K\"ahler suspensions (Proposition~\ref{fs}) and proves Theorem~\ref{mrcsplit}.
Section~\ref{furtherapplications} proves the projective bundle criterion.

\section{Disjoint Divisors and Splitting Rigidity}\label{divisors}

\subsection{Self-intersection of effective divisors}

For a smooth hypersurface, the required inequality follows from the quotient curvature formula.

\begin{proposition}\label{quoformula}
Let $D\subset X$ be a smooth complex hypersurface of a compact K\"ahler manifold $(X^n,g)$, and equip its normal bundle $N_{D/X}$ with the quotient metric.
If $e_1,\ldots,e_m$, where $m=n-1$, is a unitary tangent frame and $\nu$ is a unit complex normal, then
\[
 \sum_{a=1}^m\Theta_{N_{D/X}}(e_a,\bar e_a)
 =\Ric_g^\perp(\nu)+|\II|^2,
\]
where $\II$ is the complex second fundamental form.
\end{proposition}

In particular, if $\Ric_g^\perp\geq0$, then every smooth divisor $D$ satisfies
\[
\begin{aligned}
\int_X[D]^2\wedge\omega^{n-2}
&=\int_D c_1(N_{D/X})\wedge\omega_D^{n-2}\\
&=\frac{1}{2\pi(n-1)}
  \int_D\bigl(\Ric_g^\perp(\nu)+|\II|^2\bigr)\,\omega_D^{n-1}
\geq0.
\end{aligned}
\]
The next proposition extends this inequality to arbitrary effective divisors.

\begin{proposition}\label{Dgeq0}
Let $(X^n,g,\omega)$ be a compact K\"ahler manifold with $\Ric_g^\perp\geq0$.
Then every effective divisor $D$ satisfies
\[
\int_X[D]^2\wedge\omega^{n-2}\geq0.
\]
\end{proposition}

\begin{proof}
We first treat the case of a prime divisor $D$ and write $L_D=\cO_X(D)$.
Let $\Gamma$ be the closure of the conormal Gauss graph
\[
\Gamma^\circ
=\bigl\{(x,[du_x]):x\in D_{\mathrm{reg}}\bigr\}
\subset\Pj(T^{*1,0}X),
\]
where $u$ is a local reduced defining equation of $D$.
Choose a resolution $\mu\colon W\to\Gamma$ that is an isomorphism over $\Gamma^\circ$; see \cite{Nobile,BierstoneMilman}.
Composing $\mu$ with the natural projection to $X$ gives a proper bimeromorphic morphism
\[
p\colon W\to D\subset X
\]
that identifies $p^{-1}(D_{\mathrm{reg}})$ with $D_{\mathrm{reg}}$.
Writing $\mu(z)=(p(z),[\xi_z])$, we obtain a tautological line subbundle $K\subset p^*T^{*1,0}X$ with fibers
\[
K_z=\mathbb C\xi_z\subset(p^*T^{*1,0}X)_z,
\qquad z\in W.
\]
Its annihilator defines a holomorphic quotient line bundle
\[
Q=p^*T^{1,0}X/\operatorname{ann}(K)\simeq K^*.
\]
Equip $Q$ with the quotient metric $h_Q$ induced by $p^*g$.
Over $p^{-1}(D_{\mathrm{reg}})$, it is naturally identified with the pullback of $N_{D_{\mathrm{reg}}/X}$.

Differentiating local defining equations of $D$ gives a nonzero morphism
\[
p^*L_D^{-1}\longrightarrow K.
\]
If $E_{\mathrm{def}}\geq0$ denotes its zero divisor, then
\[
p^*L_D\simeq Q\otimes\cO_W(E_{\mathrm{def}}).
\]
Set $\beta=p^*\omega$.
On $p^{-1}(D_{\mathrm{reg}})$, the quotient curvature formula gives
\[
c_1(Q,h_Q)\wedge\beta^{n-2}
=\frac{1}{2\pi(n-1)}
  p^*\!\left(
    \bigl(\Ric_g^\perp(\nu)+\|\II_D\|^2\bigr)\omega_D^{n-1}
  \right)
\geq0.
\]
The left-hand side is a smooth top-degree form on $W$, so its nonnegativity extends to all of $W$ by density.
Since $\beta$ is semipositive, the projection formula now yields
\[
\int_X[D]^2\wedge\omega^{n-2}
=\int_W c_1(Q,h_Q)\wedge\beta^{n-2}
  +\int_{E_{\mathrm{def}}}\beta^{n-2}\\
\geq0.
\]

Now write an arbitrary nonzero effective divisor as
$D=\sum_{i=1}^k a_iZ_i$, where the $Z_i$ are distinct prime divisors and $a_i\in\mathbb Z_{>0}$.
For $i\neq j$, choose a resolution $p_i\colon W_i\to Z_i$ and pull back the canonical section of $\cO_X(Z_j)$.
Because $Z_i\not\subset Z_j$, the resulting section of $p_i^*\cO_X(Z_j)$ is not identically zero.
Writing $E_{ij}$ for its effective zero divisor, we obtain
\[
\begin{aligned}
\int_X[Z_i]\wedge[Z_j]\wedge\omega^{n-2}=\int_{W_i}p_i^*[Z_j]\wedge(p_i^*\omega)^{n-2}=\int_{E_{ij}}(p_i^*\omega)^{n-2}\geq0.
\end{aligned}
\]
Combining these inequalities with the prime-divisor case gives
\[
\begin{aligned}
\int_X[D]^2\wedge\omega^{n-2}=\sum_i a_i^2\int_X[Z_i]^2\wedge\omega^{n-2}+2\sum_{i<j}a_i a_j
  \int_X[Z_i]\wedge[Z_j]\wedge\omega^{n-2}
\geq0.
\end{aligned}
\]
\end{proof}

\subsection{A criterion for parallel kernels}

\begin{lemma}\label{parallelkernel}
Let $(X^n,g,\omega)$ be a compact K\"ahler manifold with $\Ric_g^\perp\geq0$.
Suppose that a nonzero smooth real $(1,1)$-form $\theta$ can locally be written as
\[
\theta=\sqrt{-1}\,h(q,\bar q)\,dq\wedge d\bar q,
\qquad h>0,
\]
where $q$ is holomorphic and $h$ is a smooth function of $q$ and $\bar q$.
Then the kernel of $\theta$ on
$U=\{\operatorname{tr}_\omega\theta>0\}$ extends uniquely to a parallel holomorphic hyperplane distribution $\cV$ on $X$.
\end{lemma}

\begin{proof}
The form $\theta$ is closed, semipositive, and of rank at most one.
Set $\lambda=\operatorname{tr}_\omega\theta$.
Locally, we have $\lambda=h(q,\bar q)|dq|_g^2$ and $U=\{dq\neq0\}$.
Thus $X\setminus U$ is a proper analytic subset, and $\cV=\ker\theta=\ker dq$ is a holomorphic hyperplane distribution on $U$.

Choose a local unitary frame $e_1,\ldots,e_{n-1}$ of $\cV$ and a unit normal $\nu$.
The first and mixed second derivatives of $h(q,\bar q)$ along $\cV$ vanish.
The Bochner formula for $dq$ therefore gives
\[
\begin{aligned}
\Delta_{\cV}\lambda
&=\sum_{j=1}^{n-1}
  (\partial\bar\partial\lambda)(e_j,\bar e_j)\\
&=h(q,\bar q)\sum_{j=1}^{n-1}
  (\partial\bar\partial|dq|_g^2)(e_j,\bar e_j)\\
&=h(q,\bar q)\sum_{j=1}^{n-1}
  \left(
    |\nabla_{e_j}dq|_g^2
    +|dq|_g^2R_g(e_j,\bar e_j,\nu,\bar\nu)
  \right)\\
&=h(q,\bar q)\sum_{j=1}^{n-1}|\nabla_{e_j}dq|_g^2
  +\lambda\left(
    \Ric_g(\nu,\bar\nu)-R_g(\nu,\bar\nu,\nu,\bar\nu)
  \right)\\
&=\lambda\Ric_g^\perp(\nu)
  +h(q,\bar q)\sum_{j=1}^{n-1}|\nabla_{e_j}dq|_g^2
\geq0.
\end{aligned}
\]
A direct computation also gives
\[
\sqrt{-1}\,\partial\bar\partial\lambda
\wedge\theta\wedge\omega^{n-2}
=\frac{\lambda\Delta_{\cV}\lambda}{n(n-1)}\,\omega^n
\qquad\text{on }U.
\]
The left-hand side is smooth on $X$ and vanishes on $X\setminus U$.
Since $d\theta=d\omega=0$, its integral vanishes by Stokes' theorem.
Nonnegativity then implies $\Delta_{\cV}\lambda=0$ on $U$, and the Bochner identity yields
$\nabla_u dq=0$ for every local section $u$ of $\cV$.

We claim that $\cV$ is parallel on $U$.
For a local section $u$ of $\cV$ and a $(1,0)$ vector field $Z$, the symmetry of $\nabla^{1,0}dq$ gives
\[
dq(\nabla_Zu)
=Z(dq(u))-(\nabla_Zdq)(u)
=-(\nabla_udq)(Z)
=0.
\]
Here the symmetry follows from the torsion-free property of the K\"ahler connection.
For a $(0,1)$ vector field $\bar Z$, the same conclusion follows from $\nabla_{\bar Z}dq=0$, since $dq$ is holomorphic and the $(0,1)$ part of the connection is $\bar\partial$.

It remains to extend the distribution across $X\setminus U$.
Let $P\colon T^{1,0}X|_U\to T^{1,0}X|_U$ be the orthogonal projection onto $\cV$.
The identity $\nabla P=0$ implies that $P$ is holomorphic.
Its operator norm is one, so its coefficients in local holomorphic frames are locally bounded.
Hence $P$ extends holomorphically across $X\setminus U$.
By continuity, the extension is a parallel orthogonal projection of rank $n-1$.
Its image is the required distribution, and uniqueness follows from the density of $U$.
\end{proof}

\subsection{Proof of the splitting theorem}

\begin{proof}[Proof of Theorem~\ref{main}]
\textbf{Step 1: Construct a unitary flat line bundle.}
Set
\[
\mathcal I(\xi,\zeta)
=\int_X\xi\wedge\zeta\wedge\omega^{n-2},
\qquad
t=\frac{\int_X[D_1]\wedge\omega^{n-1}}
        {\int_X[D_2]\wedge\omega^{n-1}}>0.
\]
Proposition~\ref{Dgeq0} and the disjointness of the supports give
\[
\mathcal I([D_i],[D_i])\geq0,
\qquad
\mathcal I([D_1],[D_2])=0.
\]
Since $[D_1]-t[D_2]$ is primitive, the Hodge--Riemann relations \cite[Chapter~3, \S3.3]{Huybrechts} yield
\[
0\leq\mathcal I([D_1],[D_1])+t^2\mathcal I([D_2],[D_2])\\
 =\mathcal I([D_1]-t[D_2],[D_1]-t[D_2])
 \leq0.
\]
It follows that $[D_1]=t[D_2]$.
Both divisor classes are nonzero and rational, so $t\in\mathbb Q_{>0}$.
Writing $t=a/b$ with $a,b\in\mathbb Z_{>0}$, we conclude that
\[
L=\cO_X(bD_1-aD_2)
\]
has vanishing real first Chern class and is therefore unitary flat.

\medskip
\noindent\textbf{Step 2: A semipositive form and a parallel distribution.}
Let $s$ be the canonical meromorphic section of $L$.
On a neighborhood $V$, write $s=\varphi_Ve_V$, where $e_V$ is a parallel unitary frame and $\varphi_V$ is meromorphic.
On $V\setminus\Supp(D_2)$, define
\begin{equation}\label{fs1}
\theta_V
=\frac{\sqrt{-1}\,d\varphi_V\wedge d\overline{\varphi_V}}
       {(1+|\varphi_V|^2)^2}.
\end{equation}
Because $D_1$ and $D_2$ have disjoint supports, $\varphi_V^{-1}$ is holomorphic near $\Supp(D_2)$.
Thus the local expressions
\[
\Psi_V(x)=
\begin{cases}
[1:\varphi_V(x)],
  &\text{where $\varphi_V$ is holomorphic},\\[2mm]
[\varphi_V^{-1}(x):1],
  &\text{where $\varphi_V^{-1}$ is holomorphic}
\end{cases}
\]
define a holomorphic map $\Psi_V\colon V\to\Pj^1$.
The identity $\theta_V=\Psi_V^*\omega_{\mathrm{FS}}$ shows that $\theta_V$ extends smoothly across $\Supp(D_2)$.

The transition functions between parallel unitary frames take values in $U(1)$ and preserve the Fubini--Study form.
Consequently, the forms in \eqref{fs1} glue to a nonzero smooth real $(1,1)$-form $\theta$ on $X$.
Applying Lemma~\ref{parallelkernel}, we obtain a parallel holomorphic hyperplane distribution $\cV$.
The de Rham decomposition theorem then gives the holomorphic isometric splitting
\[
\begin{gathered}
(\widetilde X,\widetilde g)
\simeq
(\Sigma,h_\Sigma)\times(N,g_N),\quad
\dim_{\mathbb C}\Sigma=1,
\quad
\widetilde{\cV}=T^{1,0}N.
\end{gathered}
\]

\medskip
\noindent\textbf{Step 3: Geometry of the reduced divisors.}
Let $\pi\colon\widetilde X\to X$ be the covering map.
Since $L$ is unitary flat, $\pi^*L$ is trivial, and the pullback of $s$ defines a holomorphic map
\[
\widetilde q\colon\widetilde X\longrightarrow\Pj^1.
\]
Its differential vanishes on $\widetilde{\cV}=T^{1,0}N$, so
$
\widetilde q=\varphi\circ\operatorname{pr}_\Sigma
$
for a nonconstant holomorphic map $\varphi\colon\Sigma\to\Pj^1$.
It follows that
\[
\pi^{-1}(\Supp D_1)=\varphi^{-1}(0)\times N, \quad
\pi^{-1}(\Supp D_2)=\varphi^{-1}(\infty)\times N.
\]
The zeros and poles of $\varphi$ are isolated, so the connected components of the lifted reduced divisors are smooth product slices $\{z_0\}\times N$.
These slices are totally geodesic, and their normal bundles have flat Chern connection with respect to the quotient metric.
These properties descend to $X$.
Finally, a parallel distribution tangent to either reduced divisor is determined by its value at a single point of that divisor, which proves uniqueness.
\end{proof}

The converse to Theorem~\ref{main} fails in general: a split universal cover does not by itself imply the existence of disjoint effective divisors.

\begin{example}\label{count1}
Let $T$ be a complex two-torus with $H^{1,1}(T)\cap H^2(T,\mathbb Z)=0$.
With a flat K\"ahler metric, it satisfies $\Ric^\perp\equiv0$ and has universal cover $\widetilde T\simeq\mathbb C\times\mathbb C$.
Nevertheless, $T$ has no nonzero effective divisors: any such divisor $D$ would have $c_1(\cO_T(D))=0$, contradicting $\int_D\omega>0$.
\end{example}

\begin{proof}[Proof of Corollary~\ref{quasipositive}]
Suppose that $X$ contains two disjoint nonzero effective divisors.
By Theorem~\ref{main}, its universal cover splits holomorphically and isometrically as $\Sigma\times N$.
At any lift $\widetilde p$ of $p$, the curvature assumption gives $\Ric_{\widetilde g}^\perp(v)>0$ for every unit vector $v$.
On the other hand, the mixed curvature terms of a product vanish, so
\[
\Ric_{\widetilde g}^\perp(v)=0
\]
for every unit vector tangent to the curve factor $\Sigma$.
This is a contradiction.
\end{proof}

\section{Applications to MRC Splitting}\label{fibrations}

We now apply Theorem~\ref{main} to MRC fibrations.
The first step is to explain how a parallel distribution with a compact rationally connected leaf determines a K\"ahler suspension.

\subsection{From parallel distributions to K\"ahler suspensions}

Let $(X,Y,F,f,\rho)$ be a locally constant fibration, and fix an identification
$
\Phi\colon (\widetilde Y\times F)/\pi_1(Y)\xrightarrow{\sim}X
$
over $Y$.
The product decomposition
\[
T^{1,0}(\widetilde Y\times F)
=\operatorname{pr}_1^*T^{1,0}\widetilde Y
 \oplus\operatorname{pr}_2^*T^{1,0}F
\]
is preserved by the diagonal action of $\pi_1(Y)$ and therefore induces a holomorphic splitting
\[
T^{1,0}X=\mathcal H\oplus\mathcal E.
\]
Since $\Phi$ commutes with the projections to $Y$, as expressed by the diagram
\[
\begin{array}{ccc}
(\widetilde Y\times F)/\pi_1(Y) & \xrightarrow{\Phi} & X\\
\downarrow & & \downarrow f\\
Y & = & Y,
\end{array}
\]
the vertical summand is $\mathcal E=T^{1,0}_{X/Y}$.
The following lemma gives a criterion for this suspension presentation to be compatible with a given K\"ahler metric.

\begin{lemma}\label{ms}
Let $(X,Y,F,f,\rho)$ be a locally constant fibration, and let $g$ be a K\"ahler metric on $X$.
Let $\mathcal H$ be the horizontal subbundle induced by the suspension presentation above, so that
\[
T^{1,0}X=\mathcal H\oplus T^{1,0}_{X/Y}.
\]
If $\mathcal H$ is orthogonal to $T^{1,0}_{X/Y}$ with respect to $g$, then there exist K\"ahler metrics $g_Y$ on $Y$ and $g_F$ on $F$ for which $(X,Y,F,f,\rho)$ is a K\"ahler suspension.
\end{lemma}

\begin{proof}
Pull back $g$ to $\widetilde Y\times F$ through the suspension presentation.
Orthogonality and the K\"ahler condition imply that the two metric blocks depend only on the coordinates of their respective factors.
Thus the pullback metric has the product form
\[
\operatorname{pr}_1^*\widetilde g_Y+\operatorname{pr}_2^*g_F.
\]
The diagonal action of $\pi_1(Y)$ preserves this metric and the product splitting.
Consequently, $\widetilde g_Y$ descends to a K\"ahler metric $g_Y$ on $Y$, and $\rho$ preserves $g_F$.
The suspension identification is therefore a holomorphic isometry.
\end{proof}

Combining Lemma~\ref{ms} with H\"oring's results on foliations with rationally connected leaves gives the following criterion.

\begin{proposition}\label{fs}
Let $(X,g)$ be a compact K\"ahler manifold with a parallel orthogonal splitting
\[
T^{1,0}X=\mathcal E\oplus\mathcal H.
\]
Suppose that the foliation defined by $\mathcal E$ has a compact rationally connected leaf.
Then there is a holomorphic submersion $f\colon X\to Y$ onto a compact K\"ahler manifold $Y$, with compact rationally connected fiber $F$, such that
$
T^{1,0}_{X/Y}=\mathcal E.
$
Moreover, there exist a K\"ahler metric $g_F$ on $F$ and a representation
$
\rho\colon\pi_1(Y)\to\Isomhol(F,g_F)
$
for which $(X,Y,F,f,\rho)$ is a K\"ahler suspension.
\end{proposition}

\begin{proof}
Since $\mathcal E$ is parallel, both $\mathcal E$ and its orthogonal complement $\mathcal H$ are integrable holomorphic subbundles.
By \cite[Corollary~2.11]{Horing}, there is a holomorphic submersion $f\colon X\to Y$ with connected fibers such that $T^{1,0}_{X/Y}=\mathcal E$.
The complementary distribution $\mathcal H$ defines an integrable holomorphic connection.
Hence \cite[Theorem~3.17]{Horing} gives a locally constant fibration $(X,Y,F,f,\rho)$ whose horizontal distribution is $\mathcal H$.
Lemma~\ref{ms} then shows that this fibration is a K\"ahler suspension.
\end{proof}

\subsection{The case \texorpdfstring{$\rd(X)=n-1$}{rd(X)=n-1}}

When the MRC base is a curve, two distinct smooth general fibers provide the disjoint divisors needed to apply Theorem~\ref{main}.

\begin{proof}[Proof of Theorem~\ref{mrcsplit}(i)]

Let $f\colon X\to C$ be the holomorphic extension of the MRC fibration.
Choose two distinct smooth general fibers $D_1$ and $D_2$.
Theorem~\ref{main} and Proposition~\ref{fs} yield a K\"ahler suspension $(X,C',F,f_0,\rho_0)$.

It remains to identify $f_0$ with $f$.
Since $C$ contains no rational curves and the fibers of $f_0$ are rationally connected,
$
f=\varphi\circ f_0
$
for a nonconstant holomorphic map $\varphi\colon C'\to C$.
The general fibers of $f$ are connected, so $\deg\varphi=1$ and $\varphi$ is a biholomorphism.
After identifying $C'$ with $C$, we obtain the required K\"ahler suspension structure.
\end{proof}

The base curve in Theorem~\ref{mrcsplit}(i) need not admit a flat metric.
The following example also illustrates that nonnegative orthogonal Ricci curvature allows negative Ricci curvature in some directions.

\begin{example}\label{nonflat}
Let $C$ be a compact hyperbolic Riemann surface, and equip $\Pj^1$ and $C$ with metrics of constant holomorphic sectional curvature $1$ and $-1$, respectively.
The product
\[
M=\Pj^1\times C
\]
then satisfies $\Ric^\perp\equiv0$.
Its MRC base is $C$, which admits no flat K\"ahler metric, and the Ricci curvature of $M$ is negative in the $C$-directions.
\end{example}

\subsection{The rationally connected case}

We next prove the equivalence in Theorem~\ref{mrcsplit}(ii).

\begin{proof}[Proof of Theorem~\ref{mrcsplit}(ii)]
Suppose that condition (a) holds.
A rationally connected compact K\"ahler manifold is simply connected, so Theorem~\ref{main} gives a holomorphic isometric product
\[
X=\Sigma\times F.
\]
Both factors are compact and simply connected. Hence $\Sigma\simeq\Pj^1$.
Moreover, $F$ is rationally connected because it is a holomorphic image of $X$.
This proves (a)$\Rightarrow$(c).

Projection onto $\Pj^1$ gives (c)$\Rightarrow$(b).
Finally, the inverse images of two distinct points under a nonconstant holomorphic map $X\to\Pj^1$ are nonzero effective divisors with disjoint supports, proving (b)$\Rightarrow$(a).
\end{proof}

\section{Further Applications}\label{furtherapplications}

\begin{proposition}\label{projectivebundle}
Let $M$ be a compact K\"ahler manifold, let $L$ be a holomorphic line bundle on $M$, and set
\[
X=\Pj(\cO_M\oplus L).
\]
Then the following conditions are equivalent:
\begin{enumerate}
\item $X$ admits a K\"ahler metric with $\Ric^\perp\geq0$.
\item $M$ admits a K\"ahler metric with $\Ric^\perp\geq0$ and $L$ is unitary flat.
\end{enumerate}
\end{proposition}

\begin{proof}
The zero and infinity sections $\sigma_0,\sigma_\infty\colon M\to X$, represented fiberwise by $[1:0]$ and $[0:1]$, have disjoint images $D_0$ and $D_\infty$.
Their normal bundles are
\[
N_{D_0/X}\simeq L,
\qquad
N_{D_\infty/X}\simeq L^{-1}.
\]

Suppose first that $X$ carries a K\"ahler metric $g$ with $\Ric_g^\perp\geq0$.
Theorem~\ref{main} implies that $L$ is unitary flat.
Since $D_0$ is totally geodesic and its normal distribution is parallel,
\[
\Ric_{g|_{D_0}}^\perp(v)=\Ric_g^\perp(v)\geq0.
\]
The identification $M\simeq D_0$ therefore gives the required metric on $M$.

Conversely, suppose that $g_M$ satisfies $\Ric_{g_M}^\perp\geq0$ and that $L$ is defined by a unitary character $\chi\colon\pi_1(M)\to U(1)$.
Then
\[
X\simeq(\widetilde M\times\Pj^1)/\pi_1(M),
\qquad
\gamma(z,[u:v])=(\gamma z,[u:\chi(\gamma)v]).
\]
Equip $\Pj^1$ with a Fubini--Study metric $h$ of constant holomorphic sectional curvature $\kappa>0$.
For unit vectors $u\in T^{1,0}\Pj^1$ and $v\in T^{1,0}M$, and for $a,b\geq0$ with $a+b=1$, the product curvature formula gives
\[ \operatorname{Ric}_{h\oplus\tilde g_M}^{\perp} (\sqrt a\,u+\sqrt b\,v) = b^2\operatorname{Ric}_{\tilde g_M}^{\perp}(v) + ab\left(\kappa+\operatorname{Ric}_{\tilde g_M}(v,\bar v)\right).  \]
Since $M$ is compact, we may choose $\kappa>0$ sufficiently large that $h\oplus\tilde g_M$ has nonnegative orthogonal Ricci curvature.
It is invariant under the above action and hence descends to the required metric on $X$.
\end{proof}

\begin{theorem}
Let $C$ be a compact Riemann surface, and let
$E\to C$ be a holomorphic vector bundle of rank $r\geq2$.
Set $X=\mathbb P(E)$.
Then the following conditions are equivalent:
\begin{enumerate}
\item $E$ is polystable.
\item $X$ admits a K\"ahler metric $g$ with nonnegative orthogonal
Ricci curvature.
\end{enumerate}
\end{theorem}

\begin{proof}
Recall that $E$ is polystable if and only if $\mathbb P(E)\to C$ admits a flat projective unitary structure.
Thus it suffices to prove that this condition is equivalent to \textup{(ii)}.

Assume \textup{(ii)}, and choose a K\"ahler form $\omega_C$ on $C$.
By Lemma~\ref{parallelkernel} applied to $\theta=p^*\omega_C$, the vertical tangent bundle $T^{1,0}_{X/C}$ is parallel.
Proposition~\ref{fs} then shows that $(X,g)$ is a K\"ahler suspension over $C$, with monodromy
\[
\rho\colon\pi_1(C)\longrightarrow
\operatorname{Isom}_{\mathrm{hol}}
(\mathbb P^{r-1},g_F).
\]
Since the holomorphic isometry group of $(\mathbb P^{r-1},g_F)$ is a compact subgroup of $\operatorname{PGL}(r,\mathbb C)$, it is conjugate to a subgroup of $\operatorname{PU}(r)$.
Thus $\mathbb P(E)$ is projectively unitary flat, and the asserted rigidity statement follows.

Conversely, suppose that $\mathbb P(E)$ is projectively unitary flat.
Then
\[
X\simeq
(\widetilde C\times\mathbb P^{r-1})/\pi_1(C),
\]
where the monodromy takes values in $\operatorname{PU}(r)$.
Choose a constant-curvature K\"ahler metric $g_C$ on $C$ with Gaussian curvature at least $-4$, and let $g_{\mathrm{FS}}$ be the Fubini--Study metric on $\mathbb P^{r-1}$ with holomorphic sectional curvature $4$.
The curvature computation in \cite[Example~1.2]{GuZhang} shows that $\widetilde g_C\oplus g_{\mathrm{FS}}$ has nonnegative orthogonal Ricci curvature.
Since the monodromy preserves $g_{\mathrm{FS}}$, this product metric descends to $X$.
Hence \textup{(ii)} holds.
\end{proof}

\medskip
\noindent\textbf{Remark.}
The assumption that the base is a curve is essential.
Indeed, Ni, Wang, and Zheng
\cite[Example~1.4]{NiWangZheng} showed that, for
$n\geq3$, the projective bundle
$
X=\mathbb P_{\mathbb P^n}
\bigl(\mathcal O_{\mathbb P^n}^{\oplus2}
\oplus\mathcal O_{\mathbb P^n}(1)\bigr)
$
admits a K\"ahler metric with
$\operatorname{Ric}^{\perp}>0$.
However, the underlying vector bundle
$
E=\mathcal O_{\mathbb P^n}^{\oplus2}
\oplus\mathcal O_{\mathbb P^n}(1)
$
is not even semistable.

\section*{Statements and Declarations}

\textbf{AI declaration.}
The research direction and main mathematical ideas of this paper were developed by the author. AI tools were primarily used to identify potential errors and gaps in proofs, and help construct counterexamples. The author takes full responsibility for the mathematical content and correctness of the paper.

\medskip
\noindent\textbf{Conflict of interest.}
The author declares no conflict of interest.

\medskip
\noindent\textbf{Funding.}
This research did not receive any specific grant from funding agencies in the public, commercial, or not-for-profit sectors.

\bibliographystyle{amsplain}
\bibliography{Ricpaper}
\end{document}